\documentclass[smallextended]{svjour3}       % onecolumn (second format)
\usepackage{amssymb,amsmath}

\newcommand{\func}[1]{\mathop{\rm #1}}
\newcommand{\dsum}{\sum}

\begin{document}

\title{A fresh geometrical look at the general S-procedure}

\author{Michel De Lara  \and Jean-Baptiste Hiriart-Urruty}

\institute{Michel De Lara\at
  CERMICS, Ecole des Ponts, Marne-la-Vall\'{e}e, France.\\
              \email{michel.delara@enpc.fr} 
           \and
           Jean-Baptiste Hiriart-Urruty \at
Institut de Math\'{e}matiques, Universit\'{e} Paul Sabatier, Toulouse,
France.
              \email{jbhu@math.univ-toulouse.fr}
}

\date{Received: date / Accepted: date}

\maketitle

\begin{abstract}
We revisit the S-procedure for general functions with ``geometrical glasses''. We
thus delineate a necessary condition, and almost a sufficient condition, to
have the S-procedure valid. Everything is expressed in terms of convexity of
augmented sets (\textit{i.e.}, via convex hulls, conical hulls) of images
built from the data functions.
\end{abstract}

\section{Introduction}

The so-called \emph{S-procedure} takes roots in Automatic Control Theory;
an excellent survey-paper on its origin and developments in that area is [$4$%
]. In the field of Optimization, the subject has also been studied
thoroughly, beginning with the quadratic data and further with general
functions. As a result, papers concerning the S-procedure abound.
Fortunately, there are from time to time survey-papers which allow to take
stock of what has been done and what needs to be done; two such examples are
[$1$] and [$6$]. With that in mind, for the convenience of the reader who is
not necessarily \textquotedblleft immersed\textquotedblright\ in the
subject, we recall in Section~\ref{The_S-procedure_for_quadratic_functions}
some of the main known results.

The keypoint of the message conveyed in our note is the following: the
essential is not the convexity of the image set of the vector-valued mapping
obtained from all the involved real-valued functions ; it is rather the
convexity of an enlarged version of this image (via operations like adding
the positive orthant $\mathbb{R}_{+}^{q}$, or taking the
conical hull). This assumption clearly is weaker than the mere convexity of
the image itself.

The S-procedure is intimately linked with the validity of a duality result
in a certain mathematical optimization problem (see a recent overview of
that in [$9$]). This was already the main motivation in \textsc{Fradkov}'s
paper~[$3$]. But this aspect is not broached here.

Our approach is essentially \emph{geometrical}; the validity of the
necessary/sufficient conditions that we develop are expressed in terms of
convexity of sets. As expected in such contexts, the main used mathematical
tool is the separation of convex sets by hyperplanes (in finite-dimensional
vector spaces). Our main results (Theorem $2$, Theorem $3$) have
similarities with some in \textsc{Fradkov}'s old paper [$3$]; they could
have been there, as much as the method as indications around some remarks
led to them. To a certain extent, our note is a revisit and an extension of
Section~$1$ in~[$3$].

\section{The S-procedure for quadratic functions}
\label{The_S-procedure_for_quadratic_functions}

We recall here some basic results on the S-procedure when only quadratic
functions are involved.

Let $Q_{0},Q_{1},\ldots,Q_{p}$ be $1+p$ real $n\times n$ symmetric matrices,
let $c_{0},c_{1},\ldots,c_{p}\in \mathbb{R}^{n}$, let $d_{0},d_{1},\ldots,d_{p}\in 
\mathbb{R}$, let $q_{i}(\cdot)$ be the associated \emph{quadratic functions}
\[
x\in \mathbb{R}^{n}\mapsto q_{i}(x)=\frac{1}{2}\left\langle
Q_{i}x,x\right\rangle +\left\langle c_{i},x\right\rangle +d_{i}. 
\]
Here and below, $\left\langle \cdot,\cdot\right\rangle $ stands for the usual
inner-product in $\mathbb{R}^{n}$.
When $c_{i}=0$ and $d_{i}=0$, one speaks of \emph{quadratic form} $q_{i}$
instead of quadratic function. When $Q_{i}=0$, one speaks of
\emph{linear} (or \emph{affine}) \emph{function}, and of \emph{linear form}
when, moreover, $d_{i}=0$.

What is called \emph{S-procedure} in Automatic Control Theory is the
relationship between%
\[
\left. \left. 
\begin{array}{c}
\left( \mathcal{I}\right) \text{ \ \ \ \ }\left( q_{i}(x)\geqslant 0\text{
for all } i=1,2,\ldots,p\right) \Rightarrow \left( q_{0}(x)\geqslant
0\right) \\ 
\text{and} \\ 
\left( \mathcal{C}\right) \text{ \ \ \ \ \ \ }\left\{ 
\begin{array}{c}
\text{There exist } \alpha _{1}\geqslant 0,\ldots,\alpha _{p}\geqslant 0\text{
such that} \\ 
q_{0}(x)-\dsum\nolimits_{i=1}^{p}\alpha _{i}q_{i}(x)\geqslant 0\text{ for
all } x\in \mathbb{R}^{n}.%
\end{array}%
\right.%
\end{array}%
\right. \right. 
\]

The implication $\left[ \left( \mathcal{C}\right) \Rightarrow \left( 
\mathcal{I}\right) \right] $ is trivial. The issue is therefore the converse
implication. We say that the \emph{S-procedure is valid} (or \emph{favorable},
or \emph{lossless}) when this converse $\left[ \left( 
\mathcal{I}\right) \Rightarrow \left( \mathcal{C}\right) \right] $ holds
true, that is to say the equivalence between the two statements $\left( 
\mathcal{I}\right) $ and $\left( \mathcal{C}\right) $. The equivalence may
be used in its negative form, \emph{i.e.} $\left[ \left( not\text{ }%
\mathcal{I}\right) \Leftrightarrow \left( not\text{ }\mathcal{C}\right) %
\right] $, whose essential content is $\left[ \left( not\text{ }\mathcal{C}%
\right) \Rightarrow \left( not\text{ }\mathcal{I}\right) \right] .$

Let us recall some important cases when the S-procedure is known to be valid
(see [$1,6$] and references therein):

- When $p=1$, provided that there exists $x_{0}$ such that $q_{1}(x_{0})>0.$

- When all the involved functions $q_{i}$ are linear forms. In that case,
this is just the \textsc{Minkowski-Farkas} lemma (in its homogeneous form).
Indeed, to have%
\[
\left\langle a_{0},x\right\rangle -\dsum\nolimits_{i=1}^{p}\alpha
_{i}\left\langle a_{i},x\right\rangle \geqslant 0\text{ for all }x\in 
\mathbb{R}^{n} 
\]%
amounts to having $a_{0}=\dsum\nolimits_{i=1}^{p}\alpha _{i}a_{i}.$

- When all the functions $q_{i}$ involved are linear functions. In that
case, this is again the \textsc{Minkowski-Farkas} lemma (non-homogeneous
form). Indeed,%
\[
\left( \left\langle a_{i},x\right\rangle -b_{i}\geqslant 0\text{ for all }%
i=1,2,\ldots,p\right) \Rightarrow \left(\left\langle
a_{0},x\right\rangle -b_{0}\geqslant 0\right) 
\]%
is equivalent to%
\[
\left\{ 
\begin{array}{c}
\text{There exist }\alpha _{1}\geqslant 0,\ldots,\alpha _{p}\geqslant 0\text{
such that} \\ 
a_{0}=\dsum\nolimits_{i=1}^{p}\alpha _{i}a_{i}\text{ and }%
b_{0}-\dsum\nolimits_{i=1}^{p}\alpha _{i}b_{i}\leqslant 0.%
\end{array}%
\right. 
\]

\section{The S-procedure for general functions}
\label{The_S-procedure_for_general_functions}

Let $f_{0},f_{1},\ldots,f_{p}:\mathbb{R}^{n}\rightarrow \mathbb{R}$ be $1+p$
(general) functions. For such a collection of functions, we mimic the
S-procedure presented for quadratic functions. The objective is to have the
equivalence between the two next assertions:%
\[
\left. 
\begin{array}{c}
\begin{array}{c}
\left( \mathcal{I}\right) \text{ \ \ \ \ \ \ }\left( f_{i}(x)\geqslant 0%
\text{ for all }i=1,2,\ldots,p\right) \Rightarrow \left( \text{ }%
f_{0}(x)\geqslant 0\right) \\ 
\text{and}%
\end{array}
\\ 
\left( \mathcal{C}\right) \text{ \ \ \ \ \ \ }\left\{ 
\begin{array}{c}
\text{There exist }\alpha _{1}\geqslant 0,\ldots,\alpha _{p}\geqslant 0\text{
such that } \\ 
f_{0}(x)-\dsum\nolimits_{i=1}^{p}\alpha _{i}f_{i}(x)\geqslant 0\text{ for
all }x\in \mathbb{R}^{n}.%
\end{array}%
\right.%
\end{array}%
\right. 
\]%
Sometimes, the expected result is written in the following \textquotedblleft
alternative theorem\textquotedblright\ form, with%
\[
\left( not\text{ }\mathcal{I}\right) \text{ }\left\{ 
\begin{array}{c}
\text{The system of inequations }\left( f_{i}(x)\geqslant 0\text{ for all }%
i=1,2,\ldots,\text{ }p\right) \\ 
\text{and } \left(f_{0}(x)<0\right) \text{ has a solution }x\in 
\mathbb{R}^{n}\text{.}%
\end{array}%
\right. 
\]

The valid S-procedure then reads: exactly one of the two statements $\left(
not\text{ }\mathcal{I}\right) $ and $\left( \mathcal{C}\right) $ is true.

\subsection{First step: when epi-convexity enters into the picture}
\label{First_step:_when_epi-convexity_enters_into_the_picture}

For real-valued functions $\varphi _{1},\varphi _{2},\ldots,\varphi _{k}$
defined on $\mathbb{R}^{n}$, we use the standard notation $\func{Im}(\varphi
_{1},\varphi _{2},\ldots,\varphi _{k})$ for the image set $\left\{ (\varphi
  _{1}(x),\varphi _{2}(x),\ldots,\varphi _{k}(x)):x\in \mathbb{R}^{n}\right\}
\subset \mathbb{R}^{k}$.

The main result in this subsection\footnote{%
From \textsc{J.-B. Hiriart-Urruty}, \emph{A remark on the general
    S-procedure.} Unpublished technical note (2020).}
is as follows.
\begin{theorem}
{Suppose that:}

{\ - There exists } $x_{0} \in\mathbb{R}^{n}${\ such that} $f_{i}(x_{0})>0$%
{\ for all} $i=1,2,\ldots,p$,

{and}

{\ - The } epi-image {of the mapping} $\left(
f_{0},-f_{1},-f_{2},\ldots,-f_{p}\right) $, that is,\\
$\func{Im}\left( f_{0},-f_{1},-f_{2},\ldots,-f_{p}\right) +
\mathbb{R}_{+}^{p+1}$ is convex.

{Then the S-procedure is valid, that is to say:} $\left( \mathcal{I}%
\right) ${\ and} $\left( \mathcal{C}\right) ${\ are equivalent.%
}
\end{theorem}

When the \emph{epi-image} of the mapping $\left(
f_{0},-f_{1},-f_{2},\ldots,-f_{p}\right) $ is convex, we say that the mapping $%
\left( f_{0},-f_{1},-f_{2},\ldots,-f_{p}\right) $ is \emph{epi-convex.}

The first assumption is common in Optimization;
it is a \textsc{Slater}-type assumption. We refer to it hereafter as
%$% \mathcal{(S)}$:
\[
\left. \left. 
\left( \mathcal{S}\right) \text{ \ \ \ \ }
\text{There exists } x_{0} \in\mathbb{R}^{n} \text{ such that }
f_{i}(x_{0})>0\text{ for all } i=1,2,\ldots,p.
\right. \right. 
\]

\begin{remark}
  %\textbf{A general remark.}
  Suppose that $\func{Im}\left(
    g_{0},g_{1},g_{2},\ldots,g_{p}\right) $ is convex. Then so is the set\\
  $\func{Im}\left(
    g_{0},-g_{1},-g_{2},\ldots,-g_{p}\right) $ %is also convex
  (as the image of the
previous set under the linear mapping $\left(
u_{0},u_{1},u_{2},\ldots,u_{p}\right) \mapsto \left(
u_{0},-u_{1},-u_{2},\ldots,-u_{p}\right) $. Hence, the set \\
$\func{Im}\left(
g_{0},-g_{1},-g_{2},\ldots,-g_{p}\right) +\mathbb{R}_{+}^{p+1}$,
sum of two convex sets, is convex.  
\end{remark}

\begin{example}
%  \textbf{Example 1}.
  Suppose that $g_{0},g_{1},g_{2},\ldots,g_{p}$ are all
convex functions. Then, the set $\func{Im}\left( g_{0},g_{1},g_{2},\ldots,g_{p}\right) $
is not necessarily convex but $\func{Im}\left(
g_{0},g_{1},g_{2},\ldots,g_{p}\right) +\mathbb{R}_{+}^{p+1}$ is
convex, as this is easily seen from the basic definition of convexity of the 
functions~$g_{i}$'s. As a result, it comes from the main theorem above that the
S-procedure is valid whenever $f_{0}$ is convex and the $%
f_{1},f_{2},\ldots,f_{p}$ are concave. We thus recover a classical result in
convex minimization (with convex inequalities).  
\end{example}

\begin{example}
  (from [$7$, Example $3.1$]). \emph{Epi-convex mapping
but with non-convex images.}
Let $q_{0}$ and $q_{1}$ be defined on $\mathbb{R}^{2}$ as follows:%
\[
q_{0}(x,y)=2x^{2}-y^{2}\text{, }q_{1}(x,y)=x+y. 
\]
Then, $\func{Im}(q_{0},q_{1})=\left\{ (u,v)\in \mathbb{R}^{2}:u\geqslant
-2v^{2}\right\} $ is not convex. However, the epi-image $\mathcal{F=}\func{Im%
}(q_{0},-q_{1})+\mathbb{R}_{+}^{2}=\mathbb{R}^{2}$ is convex.
\end{example}

\begin{example}
  Indeed a lot of effort has been made by authors to
detect (rather strong) assumptions ensuring that an image set like $\func{Im}%
\left( g_{0},g_{1},g_{2},\ldots,g_{p}\right) $ is convex, especially with
quadratic functions~$g_{i}$'s (see [$2,$ $5,7,8$]). In addition to that, it has
recently been proved that $\func{Im}(q_{1},q_{2})+\mathbb{R}^{2}$ is convex for any pair of quadratic functions $%
(q_{1},q_{2})$ ([$2$, assertion (b) in Theorem $4.19$]). The question
remains posed for a collection of three or more quadratic functions. We
conjecture that the evoked convexity result does not hold true, but we do not
have any counterexample to offer.
\end{example}

\subsection{A further step, via geometrical interpretations of $\left( 
\mathcal{I}\right) $ and $\left( \mathcal{C}\right) $}

In this subsection, we intend to provide \emph{a geometrical exact
characterization of the statement} $\left( \mathcal{I}\right) $\emph{\ and
a \textquotedblleft close to exact\textquotedblright\ geometrical
characterization of the statement} $\left( \mathcal{C}\right) $. For that
purpose, we posit:%
\begin{eqnarray*}
-\left( \text{ }\mathbb{R}_{+}^{\ast }\times \mathbb{R}_{+}^{p}\right)
&=&\left\{ (\beta _{0},\beta _{1},\ldots,\beta _{p}):\beta _{0}<0\text{ and }%
\beta _{i}\leqslant 0\text{ for all }i\right\} \\
&=&\mathcal{K}\text{ (a polyhedral convex cone in }\mathbb{R}^{p+1}\text{);}
\\
\func{Im}\left( f_{0},-f_{1},-f_{2},\ldots,-f_{p}\right) &=&\mathcal{F}\text{
(an image set in }\mathbb{R}^{p+1}\text{, from the data).}
\end{eqnarray*}

Given a set~$S$, we denote by~$\mathrm{co}S$ its convex hull, and by~$%
\mathrm{cone}S$ its convex conical hull, that is, $\left\{
\sum_{i=1}^{k}\lambda _{i}u_{i}:k\text{ positive integer, }\lambda _{i}>0%
\text{ and }u_{i}\in S\text{ for all }i\right\} $. To link the two
definitions, we clearly have that $\mathrm{cone}S=\mathbb{R}_{+}^{\ast }(%
\mathrm{co}S)=\mathrm{co}\left( \mathbb{R}_{+}^{\ast }S\right) .$
\medskip

%The main result in this subsection is as follows.
\begin{theorem}
We have the following:
\begin{eqnarray}
\left( \mathcal{I}\right) \text{ \emph{holds true}} &\Leftrightarrow &%
\mathcal{F\cap K=\emptyset }\text{ }\Leftrightarrow \text{ }\mathbb{R}%
_{+}^{\ast }\mathcal{F\cap K=\emptyset },  \label{eq:1} \\
\left( \mathcal{C}\right) \text{ \emph{holds true}} &\Rightarrow &\mathrm{%
co}\mathcal{F\cap K=\emptyset }\Leftrightarrow \text{ }\mathrm{cone}\mathcal{%
F\cap K=\emptyset },  \label{eq:2} \\
\left( \mathrm{cone}\mathcal{F\cap K=\emptyset }\text{ \emph{and} }%
\mathcal{(S)}\right) &\Leftrightarrow &\left( \mathrm{co}\mathcal{F\cap
K=\emptyset }\text{ \emph{and} }\mathcal{(S)}\right) \text{ }\mathcal{%
\Rightarrow } \left( \mathcal{C}\right) \text{ \emph{holds true}.} 
\label{eq:3}
\end{eqnarray}
\end{theorem}

In short:

- A geometrical equivalent form of $\left( \mathcal{I}\right) $ is $\mathcal{%
F\cap K=\emptyset }$ or $\mathbb{R}_{+}^{\ast }\mathcal{F\cap K=\emptyset }$.

- Provided the (slight)\textsc{\ Slater}-type assumption $(\mathcal{S})$ is
satisfied on the functions~$f_{i}$'s, a geometrical equivalent form of $\left( 
\mathcal{C}\right) $ is either\textrm{\ }$\mathrm{co}\mathcal{F\cap
K=\emptyset }$ or $\mathrm{cone}\mathcal{F\cap K=\emptyset }$.

\bigskip

\begin{proof}
  \quad % \emph{\ Proof of Theorem }$\mathit{2}$

  \textbf{-} For the first
equivalence in~\eqref{eq:1}, maybe it is easier to consider $\left( not\text{ }%
\mathcal{I}\right) $. To have $\left( not\text{ }\mathcal{I}\right) $ means
that there exists $x\in \mathbb{R}^{n}$ such that: $f_{i}(x)\geqslant 0$ for
all $i=1,\ldots,p$, and $f_{0}(x)>0$. This exactly expresses that $\mathcal{%
F\cap K\neq \emptyset }$.

The second equivalence in~\eqref{eq:1} is clear from the relation $\mathbb{R}%
_{+}^{\ast }\mathcal{K=K}$.
\medskip

\textbf{-} We intend to prove the first implication in~\eqref{eq:2}. We use the
notation $\left\langle \cdot,\cdot\right\rangle $ for the usual inner product in $%
\mathbb{R}^{p+1}=\mathbb{R\times R}^{p}$; thus $\left\langle \alpha
,z\right\rangle =\alpha _{0}z_{0}+\alpha _{1}z_{1}+\cdots+\alpha _{p}z_{p}$
whenever $\alpha =(\alpha _{0},\alpha _{1},\ldots,\alpha _{p})\in \mathbb{%
R\times R}^{p}$ and $z=(z_{0},z_{1},\ldots,z_{p})\in \mathbb{R\times R}^{p}$.

By definition of the statement $(\mathcal{C})$ itself, there exists $\alpha
=(\alpha _{0},\alpha _{1},\ldots,\alpha _{p})\in -\mathcal{K}$
(that is to say $\alpha _{0}>0$ and $\alpha _{i}$ $\geqslant 0$ for all $%
i=1,\ldots,p$) such that%
\[
\left\langle \alpha ,z\right\rangle \geqslant 0\text{ for all }%
z=(z_{0},z_{1},\ldots,z_{p})\in \mathcal{F}. 
\]
Clearly, this is equivalent to%
\begin{equation}
\left\langle \alpha ,z\right\rangle \geqslant 0\text{ for all }%
z=(z_{0},z_{1},\ldots,z_{p})\in co\mathcal{F}.
\label{eq:4}
\end{equation}

We prove by contradiction that\textrm{\ }$\mathrm{co}\mathcal{F\cap K}$ is
empty. Therefore, suppose there exists some $\beta =(\beta _{0},\beta
_{1},\ldots,\beta _{p})$ lying in $\mathrm{co}\mathcal{F\cap K}$. Then,
according to the inequality~\eqref{eq:4} just above, we get that 
\begin{equation}
\left\langle \alpha ,\beta \right\rangle =\alpha _{0}\beta
_{0}+\dsum\nolimits_{i=1}^{p}\alpha _{i}\beta _{i}\geqslant 0.  \label{eq:5}
\end{equation}
But, by definition of $\mathcal{K}$, we have $\beta _{0}<0$ and $\beta
_{i}\leqslant 0$ for all $i=1,\ldots,p$. Thus, recalling the signs of the $%
\alpha _{i}$'s, one gets at $\left\langle \alpha ,\beta \right\rangle <0$,
which contradicts~\eqref{eq:5}.

As for the equivalence in the second part of~\eqref{eq:2}, it is clear from the
following observations: $\mathrm{cone}S=\mathbb{R}_{+}^{\ast }(\mathrm{co}S)$
and $\mathbb{R}_{+}^{\ast }\mathcal{K=K}$.
\medskip

\textbf{-} We now are going to prove that $\left( \mathrm{co}\mathcal{F\cap
K=\emptyset }\text{ and }\mathcal{(S)}\right) $ $\mathcal{\Rightarrow }$ $%
\left( \mathcal{C}\right) $.

As expected in such a context, the proof is based on a separation theorem on
convex sets. Because the two convex sets $\mathrm{co}\mathcal{F}$ and $%
\mathcal{K}$ in $\mathbb{R}^{p+1}$ do not intersect, one can separate them
properly: there exists $\alpha ^{\ast }=(\alpha _{0}^{\ast },\alpha
_{1}^{\ast },\alpha _{2}^{\ast },\ldots,\alpha _{p}^{\ast })\neq 0$ in $\mathbb{%
R\times R}^{p}=\mathbb{R}^{p+1}$ such that%
\begin{eqnarray}
\sup_{b\in \mathcal{K}}\left\langle \alpha ^{\ast },b\right\rangle
&\leqslant &\inf_{z\in \mathcal{F}}\left\langle \alpha ^{\ast
},z\right\rangle =\inf_{z\in \mathrm{co}\mathcal{F}}\left\langle \alpha
^{\ast },z\right\rangle ,  \label{eq:6} \\
\inf_{b\in \mathcal{K}}\left\langle \alpha ^{\ast },b\right\rangle
&<&\sup_{z\in \mathcal{F}}\left\langle \alpha ^{\ast },z\right\rangle . 
\label{eq:7}
\end{eqnarray}

The second property~\eqref{eq:7} is useless here, due the nonemptiness of the
interior of $\mathcal{K}$.

Due to the specific structure of $\mathcal{K}$, we deduce from $(6)$ that $%
\alpha _{i}^{\ast }\geqslant 0$ for all $i=0,1,2,\ldots,p$ and, further, $%
\sup_{b\in \mathcal{K}}\left\langle \alpha ^{\ast },b\right\rangle =0$. Now,
what is in the right-hand side of $(6)$ is just $\inf_{x\in \mathbb{R}^{n}}%
\left[ \alpha _{0}^{\ast }f_{0}(x)-\dsum\nolimits_{i=1}^{p}\alpha _{i}^{\ast
}f_{i}(x)\right] $. We therefore have proved that%
\begin{equation}
\alpha _{0}^{\ast }f_{0}(x)-\dsum\nolimits_{i=1}^{p}\alpha _{i}^{\ast
}f_{i}(x)\geqslant 0\text{ for all }x\in \mathbb{R}^{n}.  \label{eq:8}
\end{equation}

We claim that $\alpha _{0}^{\ast }>0$. If not, we would have 
\[
-\dsum\nolimits_{i=1}^{p}\alpha _{i}^{\ast }f_{i}(x_{0})\geqslant 0, 
\]%
which would come into contradiction with our \textsc{Slater}-type assumption ($%
\mathcal{S}$): \emph{\ }$f_{i}(x_{0})>0$\emph{\ }for all $i=1,2,\ldots,p$,
and $\alpha _{i}^{\ast }\geqslant 0$ for all $i=1,2,\ldots,p$ (and one of them
is $>0$).

It now remains to divide $(8)$ by $\alpha _{0}^{\ast }>0$ to get at the
desired result. $\square $
  
\end{proof}

Now, we are at the point for providing a rather general geometrical
condition ensuring the validity of the S-procedure.

\begin{theorem}
  Assume the \textsc{Slater}-type condition ($\mathcal{S}$),
  and suppose there exists a set $\mathcal{Z}\subset \mathbb{R}_{+}^{p+1}$
  containing~$0$ and such that $\mathbb{R}_{+}(\mathcal{F}+\mathcal{Z})$ is a 
  convex set. Then, the S-procedure is valid, that is to say:
  ($\mathcal{I}$) implies ($\mathcal{C}$) (hence ($\mathcal{I}$)
and ($\mathcal{C}$) are equivalent).  
\end{theorem}

The set $\mathcal{Z}$ plays the role of a \textquotedblleft
convexifier\textquotedblright\ of the extended image-set $\mathbb{R}_{+}%
\mathcal{F}$. Let us see how the made assumption covers the three following
known cases:

- (The most stringent one). When the image set $\mathcal{F}$ itself is
convex, the assumed condition simply is satisfied with $\mathcal{Z=}$ $%
\left\{ 0\right\} .$

- The epi-convex case (see
\S\ref{First_step:_when_epi-convexity_enters_into_the_picture}).
Take $\mathcal{Z=}$ $\mathbb{R}_{+}^{p+1}$ to fulfill the proposed assumption. Here, instead of
considering $\mathcal{F}$ solely, one takes its so-called \textquotedblleft
upper set\textquotedblright\ $\mathcal{F}+\mathbb{R}_{+}^{p+1}$.

- The \textquotedblleft conical convex\textquotedblright\ case, \emph{i.e.}
when $\mathbb{R}_{+}\mathcal{F}$ is convex; again the considered assumption
is verified with $\mathcal{Z=}$ $\left\{ 0\right\} .$

\bigskip

\begin{proof}
\quad %\emph{Proof of Theorem 3. }

\emph{First step. }We start from the assumption \emph{(}$\mathcal{I}$%
\emph{) }in its equivalent form $\mathcal{F\cap K=\emptyset }$ (see~\eqref{eq:1}
in Theorem~$2$). We make it a bit more general by observing that $\left( 
\mathcal{F+Z}\right) $ $\mathcal{\cap }$ $\mathcal{K=\emptyset }$ for every
set $\mathcal{Z}$ contained in $\mathbb{R}_{+}^{p+1}$. This is easy to
check, as $\mathcal{Z}$ is contained in a cone placed \textquotedblleft
oppositely\textquotedblright\ to $\mathcal{K}$. We even go further by
observing that $\mathbb{R}_{+}^{\ast }\left( \mathcal{F+Z}\right) $ $%
\mathcal{\cap }$ $\mathcal{K=\emptyset }$, since $\mathcal{K}$ is a cone.
Finally, because $0\notin \mathcal{K}$, we summarize the result of this
first step in:%
\begin{equation}
\mathbb{R}_{+}\left( \mathcal{F+Z}\right) \text{ }\mathcal{\cap }\text{ }%
\mathcal{K=\emptyset }.  \label{eq:9}
\end{equation}

\emph{Second step}. Since $0\in \mathcal{Z}$, we have $\mathcal{F\subset
F+Z}$, hence $\mathcal{F\subset }$ $\mathbb{R}_{+}\left( \mathcal{F+Z}%
\right) $. By the assumed convexity of $\mathbb{R}_{+}(\mathcal{F}+\mathcal{Z%
})$, we get at 
\begin{equation}
\mathrm{co}\mathcal{F\subset }\text{ }\mathbb{R}_{+}\left( \mathcal{F+Z}%
\right) .  \label{eq:10}
\end{equation}

\emph{Final step.} We infer from~\eqref{eq:9} and~\eqref{eq:10} that $co\mathcal{F\cap }$
$\mathcal{K=\emptyset }$. It remains to apply the result~\eqref{eq:3} in Theorem~$2$
to get at the desired conclusion ($\mathcal{C}$).
  
\end{proof}

\section{Conclusion}

We have expressed all the ingredients of the general S-procedure in purely
geometrical forms, as this was initiated in the seminal paper by \textsc{%
Fradkov} ([$3$, pages $248-251$]. In doing so, we hope to have shed a fresh
new light at this kind of results, which could help to explain or to get at
new conditions for the S-procedure to be valid.

\bigskip

\textbf{References}

1.\textbf{\ }\textsc{K. Derinkuyu} and \textsc{M. C. Pinar}, \emph{On the
S-procedure and some variants.} Math. Methods Oper. Res. $64$, n$%
%TCIMACRO{\U{b0}}%
%BeginExpansion
{{}^\circ}%
%EndExpansion
1$ ($2006$), $55-77$.

2. \textsc{F. Flores-Bazan }and\textsc{\ F. Opazo, }\emph{Characterizing
the convexity of joint-range for a pair of inhomogeneous quadratic functions
and strong duality}. Minimax Theory and its Applications, Vol. $1$, n$%
%TCIMACRO{\U{b0}}%
%BeginExpansion
{{}^\circ}%
%EndExpansion
2$ ($2016$), $257-290$.

3. \textsc{A. L. Fradkov},\emph{\ Duality theorems for certain nonconvex
extremal problems.} Siberian Math. Journal $14$ ($1973$), $247-264.$

4. \textsc{S. V. Gusev} and \textsc{A. L. Likhtarnikov}, \emph{%
Kalman-Popov-Yakubovich lemma and the S-procedure: a historical survey.}
Automation and Remote Control, Vol. $67$, n$%
%TCIMACRO{\U{b0}}%
%BeginExpansion
{{}^\circ}%
%EndExpansion
11$ ($2006$), $1768-1810.$

5.\textbf{\ }\textsc{J.-B. Hiriart-Urruty} and \textsc{M. Torki}, \emph{%
Permanently going back and forth between the \textquotedblleft quadratic
world\textquotedblright\ and the \textquotedblleft convexity
world\textquotedblright\ in optimization.} J. of Applied Math. and
Optimization $45$ ($2002$), $169-184.$

6. \textsc{I. Polik} and \textsc{T. Terlaky},\emph{\ A survey of the
S-lemma}. SIAM Review, Vol. $49$, n$%
%TCIMACRO{\U{b0}}%
%BeginExpansion
{{}^\circ}%
%EndExpansion
3$ ($2007$), $371-418.$

7. \textsc{B. T. Polyack}, \emph{Convexity of quadratic transformations
and its use in control and optimization.} J. of Optimization Theory and
Applications, Vol. $99$, n$%
%TCIMACRO{\U{b0}}%
%BeginExpansion
{{}^\circ}%
%EndExpansion
3$ ($1998$), $553-583.$

8. \textsc{M. Ramana} and \textsc{A.J. Goldman}, \emph{Quadratic maps with
convex images}. Rutcor Research Report $36-94$ (October $1994$).

9. \textsc{M. Teboulle}, \emph{Nonconvex quadratic optimization: a guided
detour.} Talk in Montpellier (September $2009$), and

\emph{Hidden convexity in nonconvex quadratic optimization}. Talk at One
World Optimization Seminar (April $2020$).

\end{document}